\documentclass[12pt]{article}
\usepackage[english]{babel}
\usepackage{amsmath}
\usepackage{amssymb}
\usepackage{amsfonts}
\usepackage{amsthm}
\usepackage{graphics,graphicx}
\usepackage{color}

 \def\EE{\vbox {\hbox to 8.9pt {I\hskip-2.1pt E\hfil}}}
 \def\e{{\rm e}}
 
\begin{document}
\begin{center}
{\Large \bf A note on  the equivalence of fractional relaxation equations   to  differential equations with varying coefficients}

\vspace{1.5mm}
{\bf  Francesco Mainardi}

\textrm{Department of Physics and Astronomy,
Bologna University and INFN
\protect\\
Via Irnerio 46, I-40126 Bologna,  Italy\\
e-mail: francesco.mainardi@bo.infn.it}
\end{center}

\def\pni{\par \noindent}
\def\vsh{\smallskip\noindent}
\def\vs{\medskip\noindent}
\def\vvs{\bigskip\noindent}
\def\vvvs{\bigskip\medskip\noindent} 
\def\vsp{\smallskip\pni}  
\def\vsn{\vsh\pni}
\def\cen{\centerline}
\def\ra{\item{a)\ }} \def\rb{\item{b)\ }}   \def\rc{\item{c)\ }}
\def\eg{{\it e.g.}\ } \def\ie{{\it i.e.}\ }
\def\sg{\hbox{sign}\,}
\def\sgn{\hbox{sign}\,}
\def\sign{\hbox{sign}\,}
\def\e{\hbox{e}}
\def\exp{\hbox{exp}}
\def\ds{\displaystyle}

\vskip-0.20truecm
\centerline{\bf Abstract}

\noindent 
In this note we show how a initial value problem for a  relaxation process governed by a  differential equation of non-integer  order with a constant  coefficient may be equivalent  to that of a differential equation  of the first  order with a varying  coefficient. This equivalence is shown for the simple fractional relaxation equation that points out the relevance of the Mittag-Leffler function in fractional calculus. 
This simple argument  may lead to the equivalence of more general processes governed by  evolution equations of fractional order with constant coefficients to processes 
governed by differential equations of integer order but with varying coefficients.
Our main  motivation is to solicit the researchers  to extend this approach to other areas of applied science in order to have a more deep knowledge of certain phenomena,
both deterministic and stochastic ones, nowadays investigated with the techniques of the fractional calculus.

\vskip-0.05truecm
\noindent
{\bf Keywords:} {Caputo fractional derivatives,  Mittag--Leffler functions, anomalous relaxation.}
%
{\bf MSC:} {26A33, 33E12, 34A08, 34C26.}
\vskip -0.05truecm
 \noindent{Published on line 9 January 2018:  {\it MATHEMATICS} Vol. 6 (2018)  Paper 8, 5pp.
DOI: 10.3390/math6010008}

\setcounter{section}{0} 

\section{Introduction}
\noindent
Let  us consider the following    relaxation  equation  
\begin{equation} 
\frac{d\Psi}{dt} = - r(t)\, \Psi(t)\,,  \quad t\ge 0\,, 
\label{relaxation}
\end{equation}
subjected to the  initial condition, for the sake of simplicity,  
\begin{equation}
\Psi(0^+) =1\,,
\end{equation} 
 where  $\Psi(t)$ and  $r(t)$ are   positive functions,
 sufficiently  well-behaved for $t\ge 0$.
In Eq. \eqref{relaxation} $\Psi(t)$ denotes a non-dimensional field variable  
 and $r(t)$  the varying relaxation coefficient.
\newpage 
 \noindent
 The solution of the above initial value problem reads
\begin{equation}
\Psi(t)= \exp [-R(t)]\,, \quad R(t) = \int_0^t \!r(t')\, dt' > 0\,.
\label{solution}
\end{equation}
 It is easy to recognize from re-arranging Eq. \eqref{relaxation} that for $t\ge0$ 
   \begin{equation}
   r(t) = - \frac{\Phi(t)}{\Psi(t)}\,, \quad \Phi(t) = \Psi^{(1)}(t) = \frac{d \Psi}{dt}(t) \,.
   \end{equation}
   The solution \eqref{solution} can be derived by solving the initial value problem by separation of variables
   \begin{equation}
\int_{1}^{\Psi(t)}\!  \frac{d\Psi (t')}{\Psi(t')} = 
\int_{0}^t \! \frac{\Phi (t')}{\Psi(t')}\, dt' = -\int_0^t\!  r(t')\, dt' =  R(t).
\end{equation}
From Eq, \eqref{solution}  we also note that 
\begin{equation}
R(t) = - \log [\Psi(t)]\,.
\end{equation}
As a matter of fact,  we have shown well-known  results that will be relevant for the next Sections.  
\section{Mittag-Leffler function as solution of the fractional relaxation process}
\noindent
Let us now consider the following initial value problem for  the so-called fractional relaxation process 
\begin{equation}
\begin{cases} _*D_t^\alpha \Psi_\alpha ((t)=- \Psi_\alpha(t), & t \geq 0, \\
\Psi_\alpha(0^+)=1, 
\end{cases}
\label{fde}
\end{equation}
with $\alpha \in (0,1]$.
Above we have labeled the field variable with $\Psi_\alpha$ to point out is dependence on $\alpha$ and  the considered the  \textsl{Caputo fractional derivative}, defined as:
\begin{equation}
_*D_t^\alpha \Psi_\alpha(t) = 
\begin{cases} 
{\ds \frac{1}{\Gamma(1-\alpha)} \int_{0}^{t} \frac{\Psi_\alpha^{(1)}(t')}{(t-t')^{\alpha }}\, dt'},
 & 0  < \alpha < 1, \\ \\
{\ds \frac{d}{dt} \Psi_\alpha(t)}\,, & \alpha=1.
\end{cases}
\label{Caputo derivative}
\end{equation} 
As found in many treatises of fractional calculus, and in particular in the 2007 survey paper  by Mainardi and Gorenflo
\cite{Mainardi-Gorenflo_FCAA07} to which the interested reader is referred for details and additional references,
the solution of the fractional relaxation problem \eqref{fde}  can be obtained by using the technique of the Laplace transform in terms of the Mittag-Leffler function.
 Indeed, we get in an obvious notation
by applying the Laplace transform to Eq. \eqref{fde}
\begin{equation}
s^\alpha \widetilde \Psi_\alpha(s) - s^{\alpha-1} = - \widetilde \Psi_\alpha(s),
\quad \hbox{hence} \quad
 \widetilde \Psi_\alpha(s)=\frac{s^{\alpha-1}}{s^\alpha + 1},
 \end{equation}
 so that 
 \begin{equation}
 \Psi_\alpha(t) = E_\alpha(-\ t^\alpha) = \sum_0^\infty (-1)^n \frac{t^{\alpha n}}
 {\Gamma (\alpha n +1)}\,.
 \label{Mittag-Leffler}
 \end{equation}
 For more details on the Mittag-Leffler function we refer to the recent treatise by Gorenflo et al. \cite{GKMS_BOOK14}. Here, for readers' convenience, we report the plots of the solution \eqref{Mittag-Leffler} for some values of the parameter 
 $\alpha \in (0,1]$.    
\begin{figure} [h!]
\centering
	\includegraphics[scale=0.40]{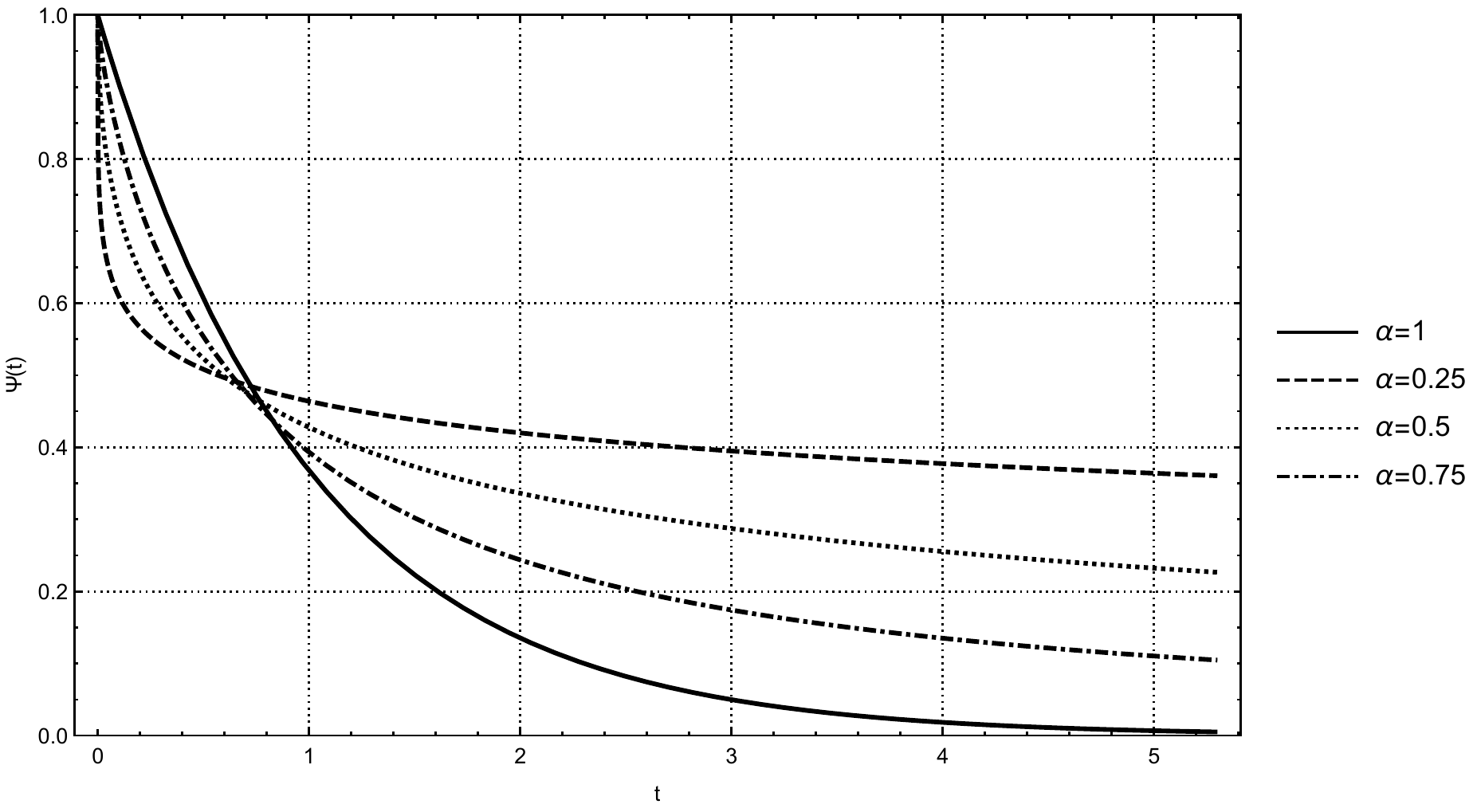}
	\vskip-0.20truecm
	\caption{Plots of the Mittag-Leffler function $\Psi_\alpha (t)$ for   
	$\alpha = 0.25, 0.50. 0.75, 1$ versus $t  \in [0.5]$.}
	\label{Fig1}
\end{figure}
It can be  noticed that for $\alpha \to 1^-$ the solution of the initial value problem reduces to the exponential function $\exp(-t)$  with a singular limit for $t \to \infty$ because
of the asymptotic representation for $\alpha \in (0,1)$,
  \begin{equation}
 E_\alpha(-\ t^\alpha) \sim \frac{t^{-\alpha}}{\Gamma(-\alpha +1)}\,, 
 \quad t \to \infty\,.
 \end{equation}
Now it is time to carry out the comparison between the two initial value problems described by Eqs. \eqref{relaxation}, \eqref{fde} with their corresponding solutions
\eqref{solution}, \eqref{Mittag-Leffler}.
It is clear that we must consider the derivative of the Mittag-Leffler function in 
\eqref{Mittag-Leffler}, namely
\begin{equation}
\Phi_\alpha(t) = \frac{d}{dt} \Psi_\alpha(t)= 
\frac{d}{dt} E_\alpha(-t^\alpha) = 
-t^{\alpha-1}\, E_{\alpha,\alpha}(-t^\alpha)\,.
\end{equation}
In Fig \ref{Fig2} we show the plots of positive function 
$-\Phi_\alpha(t)$ for some values of $\alpha \in (0,1] $.
\begin{figure} [h!]
\centering
	\includegraphics[scale=0.40]{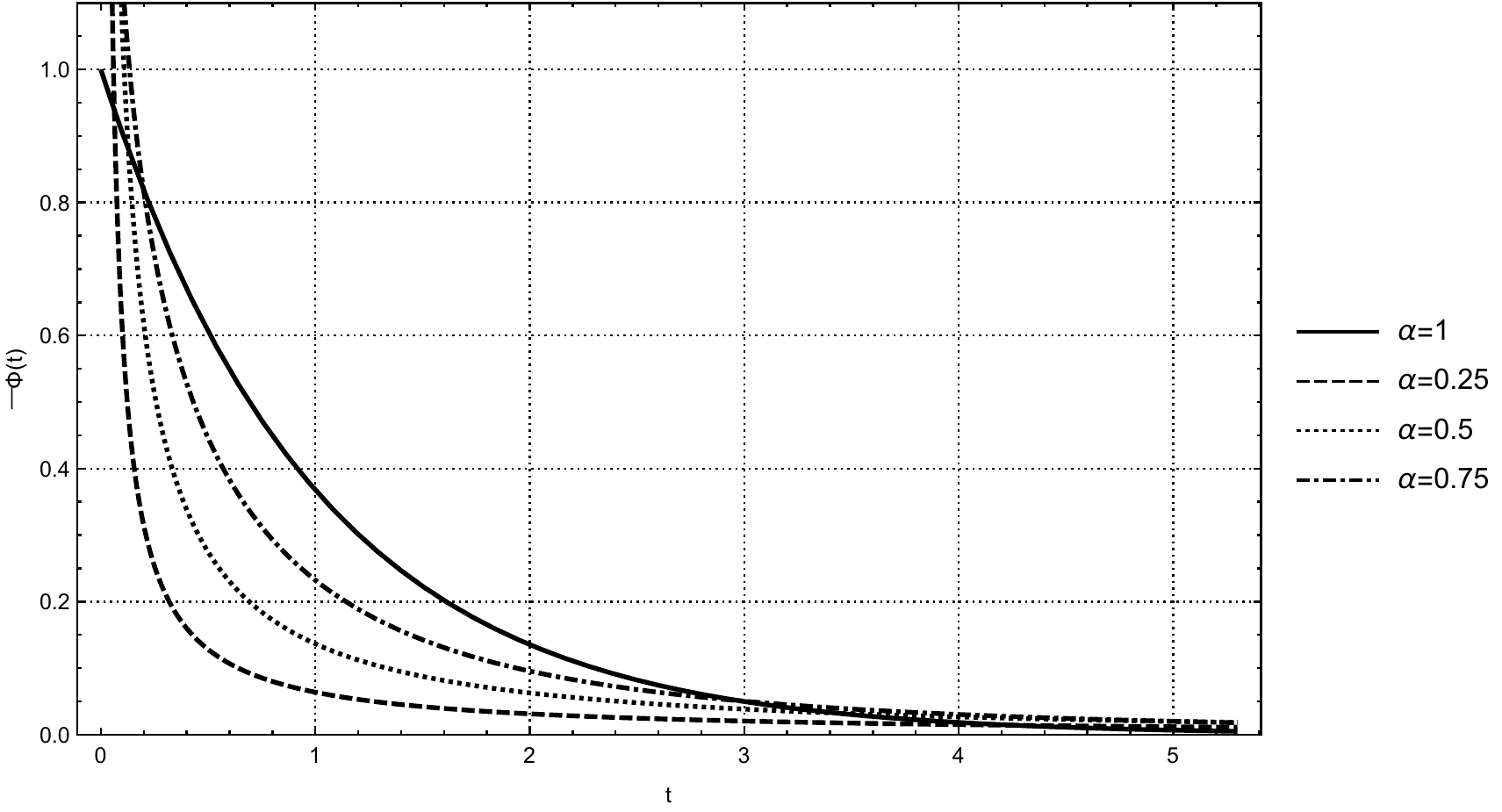}
	\vskip-0.20truecm
	\caption{Plots of the positive function  $-\Phi_\alpha (t)$ for   
	$\alpha = 0.25, 0.50. 0.75, 1$ versus $t  \in [0.5]$.}
	\label{Fig2}
\end{figure}
The above discussion leads  to the varying relaxation coefficient  of the equivalent ordinary relaxation process:  
 \begin{equation}
r_\alpha(t)= -\frac{\Phi_\alpha(t)}{\Psi_\alpha (t)} \Longrightarrow 
r_\alpha(t)= \frac{t^{\alpha -1}E_{\alpha, \alpha}(-t^\alpha)}{E_\alpha(-t^\alpha)}.
\end{equation}
Fig. \ref{r(t)} depicts the plots of $r_\alpha(t)$ for some rational values of $\alpha$, including the standard case $\alpha=1$, in which the ratio reduces to the constant 1.
\begin{figure}
\includegraphics[scale=0.40]{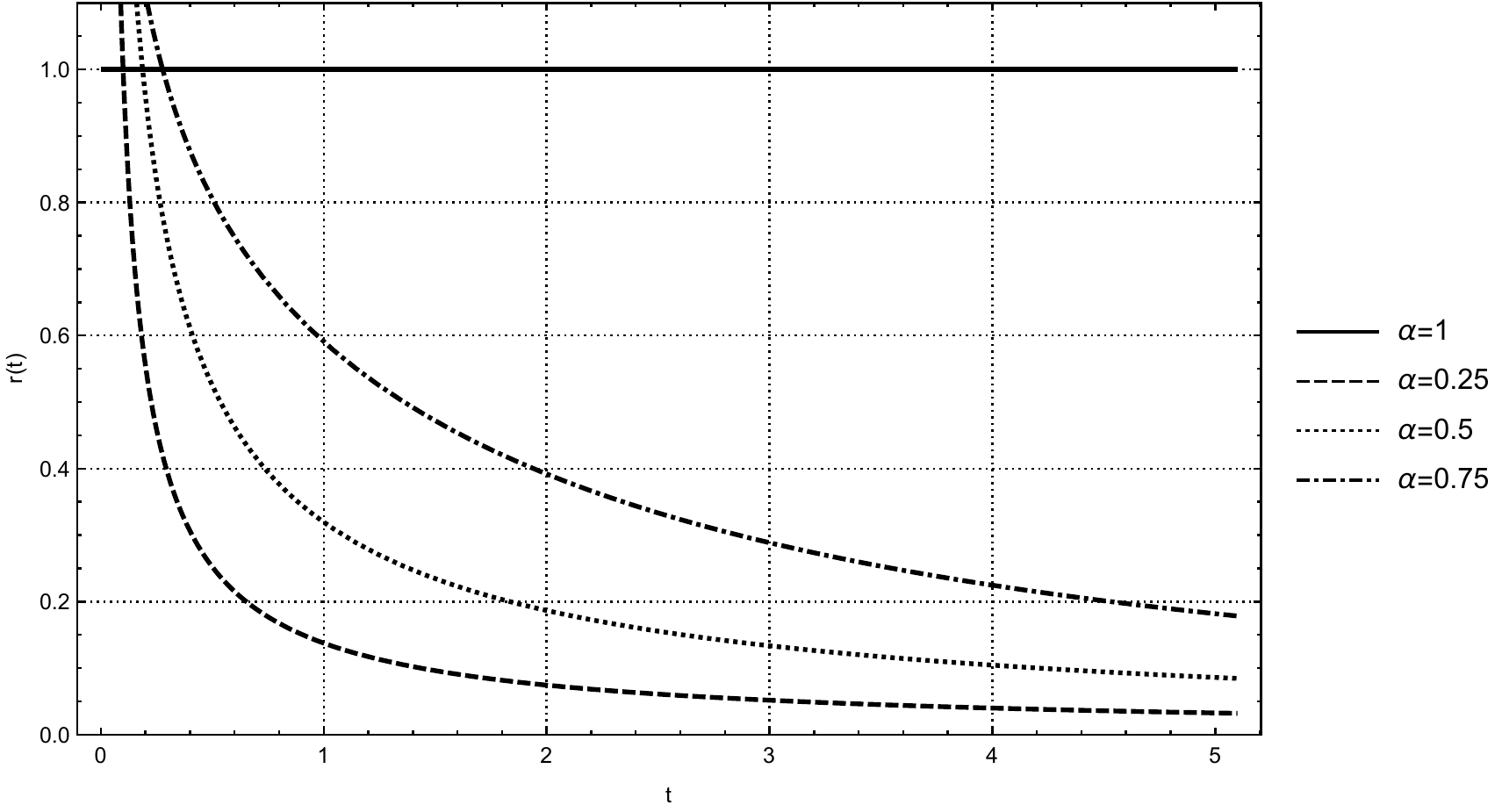}
\vskip-0.20truecm
\caption{Plots of the ratio $r_\alpha(t)$ for $\alpha =0.25, 0.50, 0.75, 1$ versus
 $t\in [0.5]$.}
\label{r(t)}
\end{figure}
\\
We conclude by plotting the function $R_\alpha(t)= -\log [\Psi_\alpha(t)]$ for some values of the parameter  $\alpha \in (0,1]$.
  \begin{figure} [h!]
\includegraphics[scale=0.40]{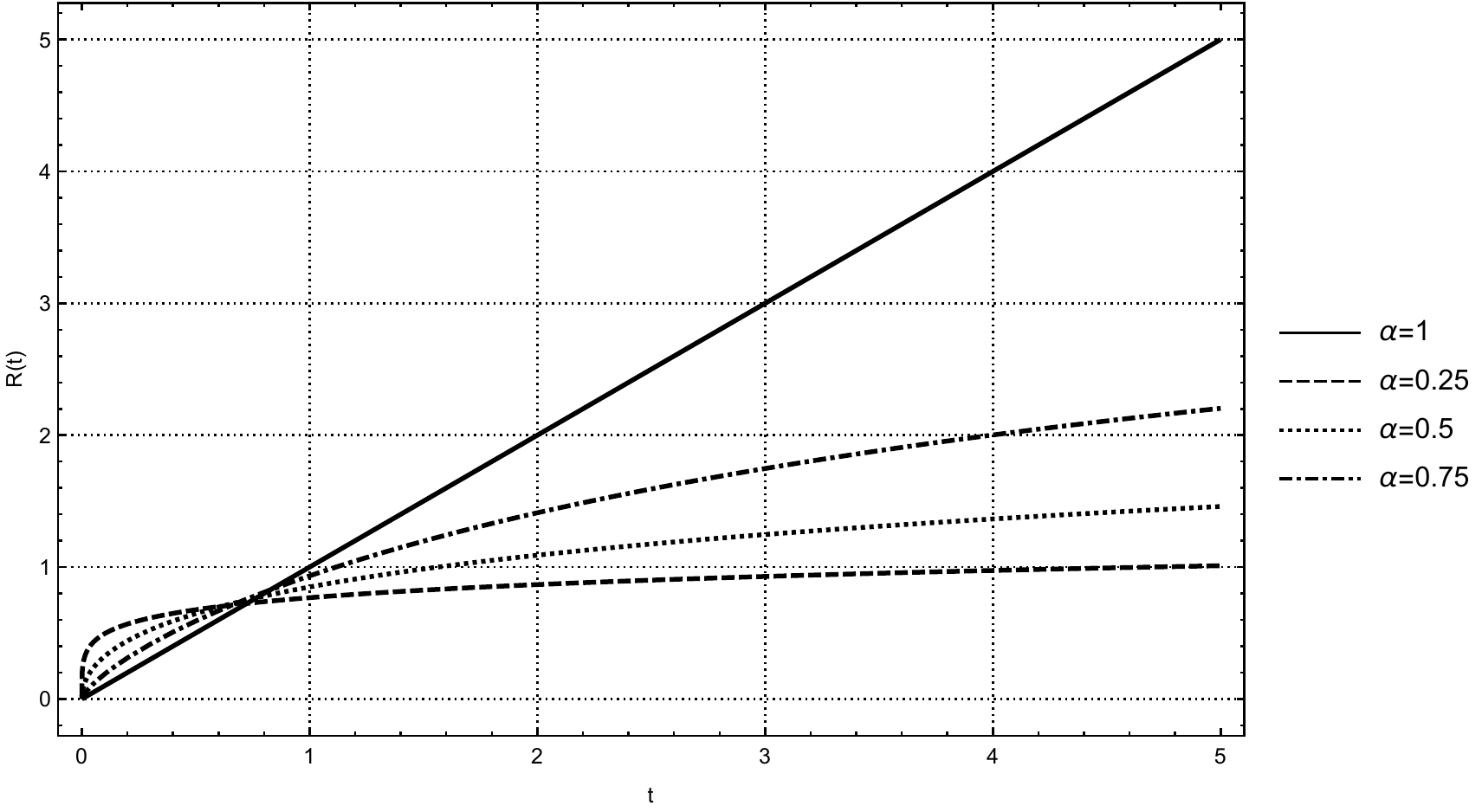}
\vskip-0.20truecm
\caption{Plots of $R_\alpha(t)$ for $\alpha=0.25, 0.50, 0.75, 1$ versus $t \in [0.5]$.}
\end{figure}

\newpage
\section{Conclusions}
\noindent
In this note we have shown how the fractional relaxation process governed by a fractional differential equation with a constant coefficient   is equivalent to a relaxation process governed by an ordinary  differential equation with a varying coefficient.
These considerations  provide a different  look at this fractional process over all for experimentalists who can measure the varying relaxation coefficient versus time.
We are convinced that it is possible  to adapt the above reasoning  to other fractional processes, including anomalous relaxation  in viscoelastic and dielectric  media and    anomalous diffusion  in complex systems. This extension is left to  
perceptive readers who  can explore these possibilities.       
\\
Last but not the least, we do not claim to be original in using the above analogy in view of the great simplicity of the argument:
for example, a similar procedure has recently been used by Sandev et al. \cite{Sandev-et-al_MATHEMATICS16} in dealing with the fractional Schr{\"o}dinger equation.    

\section*{Acknowledgments}
\noindent
The author is grateful to Leonardo Benini, student for the Master degree in  Physics
(University of Bologna),   for his valuable help in plotting.
In particular,  he has  used the MATLAB routine for the Mittag-Leffler function by Professor Roberto  Garrappa, see
\\
 https://it.mathworks.com/matlabcentral/fileexchange/48154-the-mittag-leffler-function.
 \\
 The author likes to devote this note to the memory of the late Professor Rudolf Gorenflo (1930-2017) with whom for 20 years has published joint papers. The author
 presumes that this note is written in the spirit of Prof. Gorenflo being based  on the simpler  considerations.  
\\ 
 This  work  has been carried out in the framework of the activities of the National Group of Mathematical Physics (GNFM-INdAM).

\end{document}